\documentstyle[epsfig,amsmath,12pt]{article}

\author{Petr Z\'avada\\
\\\
{\it Institute of Physics, Academy of Sciences of Czech Republic}\\
{\it Na Slovance 2, CZ-180 40 Prague 8}\\
E-mail:zavada@fzu.cz}
\title{Operator of fractional derivative in the complex plane
}
\date{July 15, 1997
}
\begin{document}

\maketitle
\begin{center}
(To be published in Communications in Mathematical Physics)
\end{center}
\begin{abstract}
The paper deals with a fractional derivative introduced by means of the
Fourier transform. The explicit form of the kernel of general derivative
operator acting on the functions analytic on a curve in complex plane is
deduced and the correspondence with some well known approaches is shown. In
particular, it is shown how the uniqueness of the operation depends on the
derivative order type (integer, rational, irrational, complex ) and the
number of poles of considered function in the complex plane.
\end{abstract}

\renewcommand{\theequation} {\thesection .\arabic{equation}}

\section{Introduction}

The fractional differentiation and integration (also called fractional
calculus) is a notion almost equally old as the ordinary differential and
integral calculus. Naturally, each a bit more gifted student who has just
understood what is the first and second derivative can ask the question: 
{\it well, but what is for example 1.5 - fold derivative?} There are several
ways how to answer such a question. An excellent review on the theory of
arbitrary order differentiation (generally of complex order) including also
interesting historical notes and comprehensive list of references to
original papers (more than thousand items) is given in the recently
published monograph \cite{samko}. Some new results were presented also in
the recent conference \cite{conf} dedicated to this topic. The fractional
calculus has also plenty of applications, see also e.g. \cite{appl},\cite
{appl1},\cite{appl2}  
and citations therein. Possible use in quantum mechanics and the
field theory is discussed in \cite{zag}. Recently, fractional derivative was
mentioned also in \cite{barci} as a particular case of pseudo - differential
operators applied in non - local field theory.

Apparently, the general prescription for a definition of fractional
derivative is using of some representation of ordinary{\it \ n-fold}
derivative (primitive function) which can be in some natural way
interpolated to{\it \ n-non\ integer.} Actually, following the mentioned
monograph, all the known approaches are always somehow connected with some
of the following relations.

1) The well known formula for {\it n-fold }integral 
\begin{equation}
\label{i1}\int_a^{x_1}dx_2\int_a^{x_2}dx_3...\int_a^{x_n}f(t)dt=\frac
1{\Gamma (n)}\int_a^{x_1}(x_1-t)^{n-1}f(t)dt 
\end{equation}
allows the substitution of $n$ by some real $\alpha >0.$ In this way
fractional integration is introduced. Then fractional derivatives can be
obtained by ordinary differentiation of fractional integrals. This is the
basis of the construction known as Riemann - Liouville fractional calculus.
Let us note, that in this approach resulting function even in the case of
fractional derivatives in general depends on the fixing of the integration
limit $a$ on the right hand side of (\ref{i1}).

2) The Cauchy formula for analytic functions in some region of the complex
plane 
\begin{equation}
\label{i2}f^{(n)}(z_0)=\frac{\Gamma (n+1)}{2i\pi }\int_C\frac{f(z)dz}{%
(z-z_0)^{n+1}} 
\end{equation}
in principle enables generalization to fractional derivatives, nevertheless
the direct extension to {\it non-integer} values of $n$ leads to
difficulties arising from multivaluedness of the term $(z-z_0)^{\alpha +1}$
and the result also depends on the choice of the cut and integration curve .

3) Analytic continuation of derivative (integral) of exponential and power
function 
\begin{equation}
\label{i3}\frac{d^\alpha }{dz^\alpha }\exp (cz)=c^\alpha \exp (cz),\qquad 
\frac{d^\alpha }{dz^\alpha }(z-c)^\beta =\frac{\Gamma (\beta +1)}{\Gamma
(\beta -\alpha +1)}(z-c)^{\beta -\alpha }. 
\end{equation}
Obviously these relations allow to define fractional derivatives of the
functions which can be expressed as linear combinations of power and
exponential functions. Also this approach is not completely consistent, as
can be illustrated by fractional derivative of exponential function expanded
to the power series 
\begin{equation}
\label{i4}\exp (cz)=\sum_{k=0}^\infty \frac{(cz)^k}{\Gamma (k+1)}, 
\end{equation}
but for $\alpha -$ {\it non integer} 
\begin{equation}
\label{i5}\frac{d^\alpha }{dz^\alpha }\exp (cz)=c^\alpha \exp (cz)\neq
c^\alpha \sum_{k=0}^\infty \frac{(cz)^{k-\alpha }}{\Gamma (k-\alpha +1)}= 
\frac{d^\alpha }{dz^\alpha }\left( \sum_{k=0}^\infty \frac{(cz)^k}{\Gamma
(k+1)}\right) . 
\end{equation}
So the task of extrapolation of integer derivative order to arbitrary one
has not the unique solution.

The approach proposed in this paper is based on the fractional derivative of
the exponential function (\ref{i3}) entering the Fourier transform of a
given function. On this basis in Sec.2 explicit form of the kernel of
fractional derivative operator is deduced. In Sec.3 the composition relation
for the derivative operator is proved. The generalization of the case of the
functions on the real axis to the case of function on the complex plane is
done in Sec.4, which is concluded by theorem summarizing the results. The
last section is devoted to the discussion of some consequences following
from the theorem and to a comparison with the known approaches as well.

\section{Definition of fractional derivative by means of Fourier transform}

\setcounter{equation}{0} Let $f(x)$ be a function having Fourier picture $
\widetilde{f}(k):$ 
\begin{equation}
\label{eq1}\widetilde{f}(k)=\int_{-\infty }^{+\infty }f(x)\exp (ikx)dx,\ \
\;\qquad f(x)=\frac 1{2\pi }\int_{-\infty }^{+\infty }\widetilde{f}(k)\exp
(-ikx)dk.
\end{equation}
Then let us create the function 
\begin{equation}
\label{eq2}f^\alpha (x)=\frac 1{2\pi }\int_{-\infty }^{+\infty }(-ik)^\alpha 
\widetilde{f}(k)\exp (-ikx)dk,\qquad \alpha >-1
\end{equation}
and define 
\begin{equation}
\label{eq3}D^\alpha (w)=\frac 1{2\pi }\int_{-\infty }^{+\infty }(-ik)^\alpha
\exp (-ikw)dk.
\end{equation}
The function $f^\alpha (x)$ can be formally expressed 
\begin{equation}
\label{eq4}f^\alpha (x)={\bf D}^\alpha f=\int_{-\infty }^{+\infty }D^\alpha
(x-y)f(y)dy.
\end{equation}
Now let us calculate the integral (\ref{eq3}), which depends on the way of
passing about the singularity $k=0$ and choice of the branch and cut
orientation of the function $k^\alpha $. For the beginning let us assume the
cut is given by the half line either $(0,-\infty )$ or $(0,+\infty )$. For
complex functions $\xi ^\alpha =(\xi _1+i\xi _2)^\alpha $ we shall accept
the phase convention 
\begin{equation}
\label{eq5}\lim _{\xi _2\rightarrow 0+}(\xi _1+i\xi _2)^\alpha =\mid \xi
_1^\alpha \mid \qquad \xi _1,\xi _2\geq 0,
\end{equation}
i.e. it holds 
\begin{equation}
\label{eq6}
\begin{array}{cccc}
& \xi <0 & \xi \geq 0 & {\rm cut\ orientation} \\ \underset{\epsilon
\rightarrow 0+}{\lim }(\xi +i\epsilon )^\alpha /\mid \xi ^\alpha \mid = & 
\left\{ 
\begin{array}{c}
\exp (+i\pi \alpha ) \\ 
\exp (+i\pi \alpha )
\end{array}
\right.  & 
\begin{array}{c}
1 \\ 
\qquad 1\qquad 
\end{array}
& 
\begin{array}{c}
(0,+\infty ) \\ 
(0,-\infty )
\end{array}
\\ 
\underset{\epsilon \rightarrow 0+}{\lim }(\xi -i\epsilon )^\alpha /\mid
\xi ^\alpha \mid = & \left\{ 
\begin{array}{c}
\exp (+i\pi \alpha ) \\ 
\exp (-i\pi \alpha )
\end{array}
\right.  & 
\begin{array}{c}
\exp (2i\pi \alpha ) \\ 
\qquad 1\qquad 
\end{array}
& 
\begin{array}{c}
(0,+\infty ) \\ 
(0,-\infty )
\end{array}
\end{array}
\end{equation}
Let us define 
\begin{equation}
\label{eq7}D_{\pm }^\alpha (w)=\frac{(-1)^\alpha }{2\pi }\int_{-\infty \pm
i0}^{+\infty \pm i0}(ik)^\alpha \exp (-ikw)dk.
\end{equation}
If we accept phase convention (\ref{eq6}) for $k^\alpha $ in the integral (%
\ref{eq7}), then the uncertainty of phase of the expression (\ref{eq7}) is
involved only in the factor 
\begin{equation}
\label{eq8}(-1)^\alpha =\exp (i\alpha \left[ 2n+1\right] \pi ),
\end{equation}
where $n$ is any integer number. Within arbitrariness given by (\ref{eq8})
the functions $D_{+}^\alpha $, $D_{-}^\alpha $ do not depend on the cut
orientation and it holds 
\begin{equation}
\label{eq9}
\begin{array}[b]{ll}
D_{+}^\alpha (w)=0 & \qquad 
{\rm for\ }w<0 \\ D_{-}^\alpha (w)=0 & \qquad {\rm for\ }w>0,
\end{array}
\end{equation}
which is evident from the corresponding integrals having the paths closed in
infinity - for $D_{+}^\alpha (D_{-}^\alpha )$ in upper (lower) half-plane.
Now we shall calculate the integrals (\ref{eq7}) in remaining regions of $w.$
Let us split them into two parts 
\begin{equation}
\label{eq10}D_{\pm }^\alpha (w)=\frac{(-1)^\alpha }{2\pi }\left[
\int_{-\infty \pm i0}^{0\pm i0}(ik)^\alpha \exp (-ikw)dk+\int_{0\pm
i0}^{+\infty \pm i0}(ik)^\alpha \exp (-ikw)dk\right] 
\end{equation}
and substitute real parameter $w$ by complex one 
\begin{equation}
\label{eq11}
\begin{array}{ll}
z_1=w+i\epsilon  & \qquad 
{\rm for\ }k<0 \\ z_2=w-i\epsilon  & \qquad {\rm for\ }k>0,
\end{array}
\end{equation}
where $\epsilon >0$. This substitution ensures absolute convergence of both
integrals in (\ref{eq10}). Next let us make in (\ref{eq10}) further
substitution 
\begin{equation}
\label{eq12}
\begin{array}{ll}
\xi =ikz_1 & \qquad 
{\rm for\ }k<0 \\ \xi =ikz_2 & \qquad {\rm for\ }k>0.
\end{array}
\end{equation}
In this way instead of $D_{\pm }^\alpha (w)$ we get functions depending also
on $\epsilon :$ 
\begin{equation}
\label{eq13}E_{\pm }^\alpha (w,\epsilon )=\frac{(-1)^\alpha }{2\pi i}\left[
\frac 1{z_1^{\alpha +1}}\int_{K_1}\xi ^\alpha \exp (-\xi )d\xi +\frac
1{z_2^{\alpha +1}}\int_{K_2}\xi ^\alpha \exp (-\xi )d\xi \right] .
\end{equation}
We assume

$a)$ Cut of complex function $k^\alpha $ is given by half-line (0,+$\infty
). $

$b)$ Function values $z_1^\alpha ,$ $z_2^\alpha $ are considered values of
one complex function $z^\alpha $ in the two different points. According to
the cut orientation of $z^\alpha $ one can get either $z_2^\alpha =($ $%
z_1^\alpha )^{*}$ or $z_2^\alpha =($ $z_1^\alpha )^{*}\exp (2i\pi \alpha )$
. Cut of the function $z^\alpha $ is assumed 
\begin{equation}
\label{eq14}
\begin{array}{ll}
(0,+\infty ) & \quad \quad 
{\rm for\ }E_{+}^\alpha  \\ (0,-\infty ) & \quad \quad {\rm for\ }%
E_{-}^\alpha .
\end{array}
\end{equation}
Latter on we shall come back to these assumptions and judge how they
affected our result.

\begin{figure}
\begin{center}
\epsfig{file=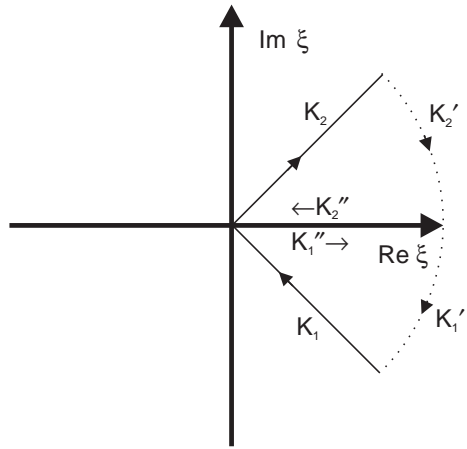,height=6cm}
\end{center}
\caption{Integration paths in Eqs. (\ref{eq13}),(\ref{eq15}).\label{fg2}}
\end{figure}

Assumptions concerning the cuts of $k^\alpha $, $z^\alpha $ correspond to
phases or to the intervals of phases of variables $k,z$ and correspondingly
to phases of the variable $\xi $. All considered possibilities are
summarized in Tab.\ref{tb1}%
\begin{table}
\begin{center}
  \begin{tabular}{|c|c|c|c|} \hline
     {\bf integral}  &  {\bf variable}  &  {\bf $k<0$}  &  {\bf $k>0$} \\ \hline
                         &         $k$            &       $\pi $         &        0      \\
        $E_{+}^\alpha $                 &          $z$             &      (0,$\pi $)    &     ($\pi $,2$\pi $)      \\
                         &           $\xi $              &  (3$\pi $/2,5$\pi $/2) & (3$\pi $/2,5$\pi $/2)  \\  \hline
                         &          $k$             &       $\pi $          &     2$\pi $         \\
   $E_{-}^\alpha $                      &          $z$             &    (0,$\pi $)         &  ($-\pi $,0)        \\
                         &        $\xi $                 & (3$\pi $/2,5$\pi $/2)    &  (3$\pi $/2,5$\pi $/2) \\ \hline
  \end{tabular}
\end{center}
   \caption{The phase correspondence of variables $\xi, z, k$ depending on the cut orientation\label{tb1}}
 \end{table}
. The corresponding integration paths are shown in Fig.\ref{fg2}. Instead of 
the interval $(3\pi /2,5\pi /2)$ for $\arg \xi $ we took the interval $(-\pi
/2,\pi /2)$ since the corresponding phase shift is the same for both
integrals in (\ref{eq13}) and can be included in factor the ($-1$)$^\alpha $
ahead of the integrals. Obviously, it holds 
\begin{equation}
\label{eq15}\int_{K_2}\xi ^\alpha \exp (-\xi )d\xi =-\int_{K_1}\xi ^\alpha
\exp (-\xi )d\xi =\int\nolimits_0^\infty \xi ^\alpha \exp (-\xi )d\xi
=\Gamma (\alpha +1).
\end{equation}
After inserting into (\ref{eq13}) we get 
\begin{equation}
\label{eq16}E_{\pm }^\alpha (w,\epsilon )=\frac{(-1)^{\alpha +1}\Gamma
(\alpha +1)}{2i\pi }\left( \frac 1{(w+i\epsilon )^{\alpha +1}}-\frac
1{(w-i\epsilon )^{\alpha +1}}\right) ,\qquad \qquad \epsilon >0.
\end{equation}
We shall regard the integrals (\ref{eq7}) as the limits 
\begin{equation}
\label{eq17}D_{\pm }^\alpha (w)=\lim _{\epsilon \rightarrow 0+}E_{\pm
}^\alpha (w,\epsilon ),
\end{equation}
hence the kernel of the operator (\ref{eq3}) which can be also considered
the generalized function is symbolically written 
\begin{equation}
\label{eq18}D_{\pm }^\alpha (w)=\frac{(-1)^{\alpha +1}\Gamma (\alpha +1)}{%
2i\pi }\left( \frac 1{(w+i0)^{\alpha +1}}-\frac 1{(w-i0)^{\alpha +1}}\right)
,
\end{equation}
where the two modes correspond to different cuts of $(w+i\epsilon )^{\alpha
+1}$ 
\begin{equation}
\label{eq19}
\begin{array}{ll}
D_{+}^\alpha (w) & \qquad \qquad 
{\rm for\ cut\ }(0,+\infty ) \\ D_{-}^\alpha (w) & \qquad \qquad {\rm for\
cut\ }(0,-\infty ).
\end{array}
\end{equation}
Let us note the function (\ref{eq18}) is well defined for any complex $%
\alpha =\alpha _1+i\alpha _2\neq -1,-2,..$, but for the beginning we assume $%
\alpha _2=0$. Now let us go back to the assumptions $a)$, $b)$ which imply
result (\ref{eq16}).

$a$)Assuming the opposite cut orientation (0,$-\infty )$ for $k^\alpha $ in
Eq. (\ref{eq10}) and repeating the corresponding sequence of steps does not
change anything for $E_{+}^\alpha (w,\epsilon )\ $whereas in the case of $%
E_{-}^\alpha (w,\epsilon )$ the phase $\arg \xi $ will shift by $-$2$\pi $
in both integrals in (\ref{eq13}), but this change can be included in factor
($-1$)$^\alpha $ ahead of the integral.

$b$) Let us assume the cuts in $z^\alpha $ having the opposite orientation
than in (\ref{eq14}). Let us take e.g. function $E_{+}^\alpha (w,\epsilon )$%
, then the phase of $\xi $ complies with 
\begin{equation}
\label{eq20}
\begin{array}{llll}
3\pi /2 & <\arg \xi < & 5\pi /2 & \qquad \quad \quad 
{\rm for\ }k>0 \\ -\pi /2 & <\arg \xi < & \pi /2 & \qquad \quad \quad {\rm %
for\ }k<0,
\end{array}
\end{equation}
i.e. the phase of second integral in (\ref{eq13}) is now shifted by $-$2$\pi
\alpha $ and instead of (\ref{eq16}) we get 
\begin{equation}
\label{eq21}E_{+}^\alpha (w,\epsilon )=\frac{(-1)^{\alpha +1}\Gamma (\alpha
+1)}{2i\pi }\left( \frac 1{(w+i\epsilon )^{\alpha +1}}-\frac{\exp (-2i\pi
\alpha )}{(w-i\epsilon )^{\alpha +1}}\right) ,\qquad \qquad \epsilon >0.
\end{equation}
Obviously this function for $w\in (-\infty ,+\infty )$ with assumed cut of $%
z^\alpha $ on $(0,-\infty )$ is identically equal to the function (\ref{eq16}%
) with the cut on $(0,+\infty )$. The analogous result can be obtained also
for $E_{-}^\alpha .$

So we have shown, that result (\ref{eq16}) of the integration (\ref{eq7})
does not depend on the choice of cut orientation of the functions $k^\alpha $
and $w^\alpha .$ The cut orientation of $w^\alpha $ in (\ref{eq16}) and (\ref
{eq18}) is dictated only by the way of passing around the singularity $k=0$
in initial integrals (\ref{eq7}) and the correspondence (\ref{eq19}) always
holds. Next let us notice the important property of the functions $D_{\pm
}^\alpha ,$ $E_{\pm }^\alpha $. Obviously it holds 
\begin{equation}
\label{eq22} 
\begin{array}{lll}
\frac{
\begin{array}{c}
d 
\end{array}
}{
\begin{array}{c}
dw 
\end{array}
}E_{\pm }^\alpha (w,\epsilon ) & = & E_{\pm }^{\alpha +1}(w,\epsilon ), \\ 
\frac{
\begin{array}{c}
d 
\end{array}
}{
\begin{array}{c}
dw 
\end{array}
}D_{\pm }^\alpha (w) & = & D_{\pm }^{\alpha +1}(w). 
\end{array}
\end{equation}
Now let us go back to the Eq. (\ref{eq4}) and consider how to calculate this
integral. It is possible either

$a)$ First to make the difference of both terms in (\ref{eq16}), then
integration and limit for $\epsilon \rightarrow 0,$

\noindent
or

$b)$ First to calculate both integrals independently, then make the limit of
their difference.

\noindent
Let us discuss the both ways and compare the results.

{\it $a)$ The action of operator ${\bf D}_{\pm }^\alpha \ $as the integral
of difference}

\noindent
We make calculation separately for the two cases:

{\it $a1)$ $\alpha =n\geq 0$ is an integer number}

\noindent
In that case the complex function $w^\alpha $ has no cuts (we can omit the
subscript $\pm )$ and for $n=0$ we can write 
\begin{equation}
\label{eq23}E^0(w,\epsilon )=\frac 1\pi \cdot \frac \epsilon {w^2+\epsilon
^2},
\end{equation}
i.e. for $\epsilon \rightarrow 0$ we get the known representation of $\delta 
$-function (see e.g. \cite{delta},p.35) 
\begin{equation}
\label{eq24}\delta (w)=\lim _{\epsilon \rightarrow 0+}\frac 1\pi \cdot \frac
\epsilon {w^2+\epsilon ^2}
\end{equation}
which acts on the function $f,$ 
\begin{equation}
\label{eq24+}\int_{-\infty }^{+\infty }\delta (x-y)f(y)dy=\lim _{\epsilon
\rightarrow 0+}\frac 1\pi \int_{-\infty }^{+\infty }\frac \epsilon
{(x-y)^2+\epsilon ^2}\ f(y)dy=f(x).
\end{equation}
Using equations (\ref{eq22})$-$(\ref{eq24+}) one can easily show 
\begin{equation}
\label{eq25}\frac{d^nf(x)}{dx^n}=\lim _{\epsilon \rightarrow
0+}\int_{-\infty }^{+\infty }E^n(x-y,\epsilon )f(y)dy
\end{equation}
for any $n\geq 0$ for which the integral converges. So it is possible to
identify 
\begin{equation}
\label{eq26}
\begin{array}{l}
D^n(w)=
\frac{
\begin{array}{c}
d^n
\end{array}
}{
\begin{array}{c}
dw^n
\end{array}
}\delta (w) \\ \frac{
\begin{array}{c}
d^nf(x)
\end{array}
}{
\begin{array}{c}
dx^n
\end{array}
}=\int_{-\infty }^{+\infty }D^n(x-y)f(y)dy,
\end{array}
\end{equation}
i.e. action of the operator ${\bf D}^n$ corresponds to {\it n-fold}
derivative.

$a2)${\it \ $\alpha $ is a real, non integer number}

\noindent
Taking into account the cut orientations (\ref{eq19}) calculation of limits (%
\ref{eq17}) gives 
\begin{equation}
\label{eq27} 
\begin{array}{ll}
D_{-}^\alpha (w)=0 & \qquad \quad w>0 \\ 
D_{+}^\alpha (w)=0 & \qquad \quad w<0 \\ 
D_{-}^\alpha (w)=-\frac{
\begin{array}{c}
\Gamma (\alpha +1)\sin ([\alpha +1]\pi ) 
\end{array}
}{
\begin{array}{c}
\pi w^{\alpha +1} 
\end{array}
} & \qquad \quad w<0 \\ 
D_{+}^\alpha (w)=+\frac{
\begin{array}{c}
\Gamma (\alpha +1)\sin ([\alpha +1]\pi ) 
\end{array}
}{
\begin{array}{c}
\pi w^{\alpha +1} 
\end{array}
} & \qquad \quad w>0. 
\end{array}
\end{equation}
The first two equations consist with Eqs. (\ref{eq9}), therefore again
confirm correct correspondence of cuts in (\ref{eq19}). After inserting of (%
\ref{eq27}) into Eq. (\ref{eq4}) we get 
\begin{equation}
\label{eq28}f_{\pm }^\alpha (x)=-\frac{\Gamma (\alpha +1)\sin ([\alpha
+1]\pi )}\pi \int_x^{\mp \infty }\frac{f(y)dy}{(x-y)^{\alpha +1}}. 
\end{equation}
For $y=x$ and $\alpha \geq 0$ the integral has a singularity. The method for
regularization of this integral will become apparent in the next.

$b)${\it The action of operator }${\bf D}${\it $_{\pm }^\alpha \ $as a
difference of two integrals}

\noindent
After inserting of (\ref{eq18}) to (\ref{eq4}) we get 
\begin{equation}
\label{eq29}f_{\pm }^\alpha (x)=\lim _{\epsilon \rightarrow 0+}\left[ \frac{%
(-1)^{\alpha +1}\Gamma (\alpha +1)}{2i\pi }\left( \int_{-\infty }^{+\infty } 
\frac{f(y)dy}{(x-y+i\epsilon )^{\alpha +1}}-\int_{-\infty }^{+\infty }\frac{%
f(y)dy}{(x-y-i\epsilon )^{\alpha +1}}\right) \right] . 
\end{equation}
We assume for the present the function $f(y)$ is analytic on the whole real
axis. Last equation can be rewritten 
\begin{equation}
\label{eq30}f_{\pm }^\alpha (x)=\lim _{\epsilon \rightarrow 0+}\left[ \frac{%
(-1)^{\alpha +1}\Gamma (\alpha +1)}{2i\pi }\left( \int_{-\infty -i\epsilon
}^{+\infty -i\epsilon }\frac{f(z)dz}{(x-z)^{\alpha +1}}-\int_{-\infty
+i\epsilon }^{+\infty +i\epsilon }\frac{f(z)dz}{(x-z)^{\alpha +1}}\right)
\right] . 
\end{equation}
Again let us separate two cases:

{\it $b1)$ $\alpha =n\geq 0$ is an integer number}

\noindent
In Eq. (\ref{eq30}) we can link up both integration paths in $z=\pm \infty $
and write 
\begin{equation}
\label{eq31}f_{}^n(x)=\lim _{\epsilon \rightarrow 0+}\frac{(-1)^{n+1}\Gamma
(n+1)}{2i\pi }\int_C\frac{f(z)dz}{(x-z)^{n+1}}=\frac{d^nf}{dx^n}
\end{equation}
where $C$ is any closed curve enclosing the singular point. Therefore for $%
n\geq 0$ the operator ${\bf D}{\it ^n}$ can be again identified with
ordinary {\it n-fold} derivative.

$b2)${\it \ $\alpha $ is a real, non integer number}

\begin{figure}
\begin{center}
\epsfig{file=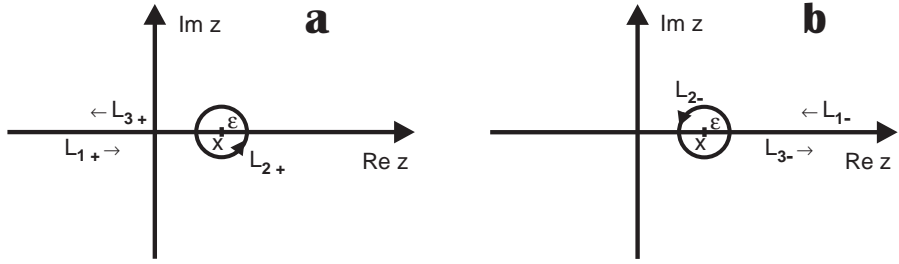,width=12cm}
\end{center}
\caption{Integration paths in Eq. (\ref{eq33}): $a)C_{+},\ b)C_{-}$.\label{fg5}}
\end{figure}

\noindent
Similarly, as in the case $b1)$ we get the integrals 
\begin{equation}
\label{eq33}f_{\pm }^\alpha (x)=\lim _{\epsilon \rightarrow 0+}\frac{%
(-1)^{\alpha +1}\Gamma (\alpha +1)}{2i\pi }\int_{C_{\pm }}\frac{f(z)dz}{%
(x-z)^{\alpha +1}}
\end{equation}
where the paths $C_{\pm }$ can pass about the corresponding cut as shown in
Fig.\ref{fg5}. The integrals converge provided that%
\begin{equation}
\label{eq33+}\lim _{z\rightarrow \pm \infty }\frac{f(z)}{z^\alpha }=0.
\end{equation}
The integration paths can be decomposed into three parts $L_{1\pm }$, $%
L_{2\pm }$, $L_{3\pm }$ and after evaluation of the corresponding integrals
one gets 
\begin{equation}
\label{eq38}f_{\pm }^\alpha (x)=-\lim _{\epsilon \rightarrow 0+}\frac{\Gamma
(\alpha +1)\sin ([\alpha +1]\pi )}\pi \left( \int_{x\mp \epsilon }^{\mp
\infty }\frac{f(z)dz}{(x-z)^{\alpha +1}}+\frac{(\pm 1)^\alpha f(x)}{\alpha
\epsilon ^\alpha }\right) ,
\end{equation}
or using the relation (\ref{A7}) from Appendix, 
\begin{equation}
\label{eq38+}f_{\pm }^\alpha (x)=-\lim _{\epsilon \rightarrow 0+}\frac
1{\Gamma (-\alpha )}\left( \int_{x\mp \epsilon }^{\mp \infty }\frac{f(z)dz}{%
(x-z)^{\alpha +1}}+\frac{(\pm 1)^\alpha f(x)}{\alpha \epsilon ^\alpha }%
\right) .
\end{equation}
For $\alpha <0$ the integrals are finite and the second terms vanish. In
this case Eq. (\ref{eq38}) is identical with Eq. (\ref{eq28}). On the other
hand we have assumed the initial integral (\ref{eq33}) is finite, therefore
the sum of both terms in (\ref{eq38}) is also finite even for any $\alpha >0.
$ In this sense by addition of the second term the integral in (\ref{eq28})
can be regularized. Note, that in contradistinction to (\ref{eq18}) the
relation (\ref{eq38+}) is well defined also for $\alpha =-1,-2,-3,...$ but $%
\alpha \neq 0,1,2,3,...$ The last equation for $\alpha $ {\it negative
integer} 
\begin{equation}
\label{eq39-}f_{\pm }^{-n}(x)=\frac 1{\Gamma (n)}\int_{\mp \infty
}^x(x-z)^{n-1}f(z)dz
\end{equation}
is a special case of the {\it n-fold} integral formula (\ref{i1}).

Now, let us assume $0<\alpha <1$ and calculate integral (\ref{eq33}) by
parts. Obviously 
\begin{equation}
\label{eq39}\int_{C_{\pm }}^{}\frac{f(z)dz}{(x-z)^{\alpha +1}}=-\frac
1\alpha \int_{C_{\pm }}\frac{f^{\,\prime }(z)dz}{(x-z)^\alpha }+\left[ \frac{%
f(z)}{(x-z)^\alpha }\right] _{\mp \infty \mp \,i0}^{\mp \infty \pm \,i0}
\end{equation}
and integral on right side is finite. Any $\alpha $ can be written as the
sum of integer and fractional part 
\begin{equation}
\label{eq40}\alpha =n+\Delta \alpha ,\qquad n=[\alpha ],\qquad 0\leq \Delta
\alpha <1.
\end{equation}
For $n\geq 0$ we can repeat integration by parts $n+1$ times and if the
function $f$ meets requirements 
\begin{equation}
\label{eq41}\lim _{z\rightarrow \pm \infty }\frac{f^{\,p}(z)}{(x-z)^{\alpha
-p}}=0\qquad {\rm for\ \ }p=0,1,...n,
\end{equation}
then instead of (\ref{eq38}) we get 
$$
f_{\pm }^\alpha (x)=\lim _{\epsilon \rightarrow 0+}\frac{(-1)^{\Delta \alpha
}\Gamma (\Delta \alpha )}{2i\pi }\int_{C_{\pm }}\frac{f^{\,n+1}(z)dz}{%
(x-z)^{\Delta \alpha }}= 
$$
\begin{equation}
\label{eq42}=-\frac{\Gamma (\Delta \alpha )\sin (\Delta \alpha \pi )}\pi
\int_x^{\mp \infty }\frac{f^{\,n+1}(z)dz}{(x-z)^{\Delta \alpha }}.
\end{equation}
Let us note that integrals (\ref{eq33}) can be modified also in other way.
We assume that all integrals 
\begin{equation}
\label{eq43}I_{\pm }(\gamma ,x)=\int_{C_{\pm }}\frac{f^{\,}(z)dz}{%
(x-z)^\gamma },\qquad \qquad \gamma =\Delta \alpha ,\ \Delta \alpha
+1,...\alpha +1
\end{equation}
converge, then recurrent relation holds 
\begin{equation}
\label{eq44}\frac d{dx}I_{\pm }(\gamma ,x)=-\gamma I_{\pm }(\gamma +1,x).
\end{equation}
Application of this relation for integral (\ref{eq33}) and $\alpha =n+\Delta
\alpha $ gives 
$$
f_{\pm }^\alpha (x)=\lim _{\epsilon \rightarrow 0+}\frac{(-1)^{\Delta \alpha
}\Gamma (\Delta \alpha )}{2i\pi }\left( \frac d{dx}\right)
^{n+1}\int_{C_{\pm }}\frac{f^{\,}(z)dz}{(x-z)^{\Delta \alpha }}= 
$$
\begin{equation}
\label{eq45}=-\frac{\Gamma (\Delta \alpha )\sin (\Delta \alpha \pi )}\pi
\left( \frac d{dx}\right) ^{n+1}\int_x^{\mp \infty }\frac{f^{\,}(z)dz}{%
(x-z)^{\Delta \alpha }}.
\end{equation}
Actually in Eq. (\ref{eq42}) we apply operation ${\bf D}_{\pm }^{\Delta
\alpha -1}$ on {\it n+1 - fold} derivative of the function $f$ whereas in
Eq. (\ref{eq45}) the both operations are interchanged. Using relation (\ref
{A7}) obviously we can modify both equations: 
\begin{equation}
\label{eq46}f_{\pm }^\alpha (x)=\frac 1{\Gamma (1-\Delta \alpha )}\int_{\mp
\infty }^x\frac{f^{\,n+1}(z)dz}{(x-z)^{\Delta \alpha }}=\frac 1{\Gamma
(1-\Delta \alpha )}\left( \frac d{dx}\right) ^{n+1}\int_{\mp \infty }^x\frac{%
f^{\,}(z)dz}{(x-z)^{\Delta \alpha }}.
\end{equation}

So we have shown the operator ${\bf D}^\alpha $ with the kernel (\ref{eq18})
can be considered a continuous interpolation of the ordinary {\it n-fold}
derivative (integral) of the functions analytic on real axis and fulfilling
the condition (\ref{eq33+}).

\section{Composition of fractional derivatives}

\begin{figure}
\begin{center}
\epsfig{file=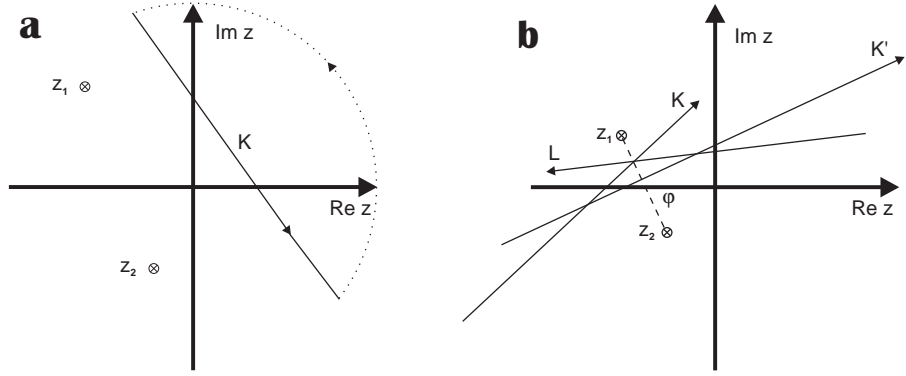,height=5cm}
\end{center}
\caption{Integration paths in Eq. (\ref{eq50}): $a)$case when the integral
vanishes, $b)$case leading to the result (\ref{eq54}).\label{fg6}} 
\end{figure}

\setcounter{equation}{0} In this section we shall investigate how the
composition of the operators {\bf D}$^\alpha $ is realized in the
representation given by Eq. (\ref{eq18}). Therefore we shall deal with the
integrals 
\begin{equation}
\label{eq47}I=\int_{-\infty }^{+\infty }D^\alpha (x-y)D^\beta (y-z)dy.
\end{equation}
Let us denote 
\begin{equation}
\label{eq48}
\begin{array}{cc}
\begin{array}{l}
h^{\bullet }(\gamma ,w)=
\frac{
\begin{array}{c}
1
\end{array}
}{
\begin{array}{c}
(w+i\tau )^\gamma 
\end{array}
} \\ h_{\bullet }(\gamma ,w)=-\frac{
\begin{array}{c}
1
\end{array}
}{
\begin{array}{c}
(w-i\tau )^\gamma 
\end{array}
}
\end{array}
& \qquad \qquad \tau >0
\end{array}
\end{equation}
and 
\begin{equation}
\label{eq48+}
\begin{array}{l}
I_1=\int_{-\infty }^{+\infty }h^{\bullet }(\alpha +1,x-y)h^{\bullet }(\beta
+1,y-z)dy \\ 
I_2=\int_{-\infty }^{+\infty }h_{\bullet }(\alpha +1,x-y)h_{\bullet }(\beta
+1,y-z)dy \\ 
I_3=\int_{-\infty }^{+\infty }h_{\bullet }(\alpha +1,x-y)h^{\bullet }(\beta
+1,y-z)dy \\ 
I_4=\int_{-\infty }^{+\infty }h^{\bullet }(\alpha +1,x-y)h_{\bullet }(\beta
+1,y-z)dy,
\end{array}
\end{equation}
then 
\begin{equation}
\label{eq49}I=\frac{(-1)^{\alpha +\beta +1}\Gamma (\alpha +1)\Gamma (\beta
+1)}{4\pi ^2}(I_1+I_2+I_3+I_4).
\end{equation}
Now let us calculate the more general integral 
\begin{equation}
\label{eq50}J=\int_K\frac{dz}{(z_2-z)^{\alpha +1}(z-z_1)^{\beta +1}},
\end{equation}
where $K$ is arbitrary line in complex plane, $z_1,z_2$ any two (diverse)
points and let $\alpha +\beta >-1.$ We also assume the cuts of the function
in the integral do not intersect the line $K.$ Then there are two
possibilities:

$a)$ $K$ is not passing between the points $z_1,z_2,$ see Fig.\ref{fg6}$a$. 
Then obviously $J=0,$ since the line $K$ can be closed in infinity by the
arc in half plane which do not contain singularities $z_1,z_2.$

$b)$ $K$ is passing between the points $z_1,z_2,$ see Fig.\ref{fg6}$b$. $%
K^{\prime }$ denotes line crossing the segment $\left\langle
z_1,z_2\right\rangle $ perpendicularly at its center. If we assume, that the
cuts do not intersect any of both lines $K,K^{\prime },$ then in the
integral (\ref{eq50}) path $K$ can be substituted by $K^{\prime }.$ Further,
if we denote 
\begin{equation}
\label{eq51}z_0=(z_1+z_2)/2,\qquad \qquad r\exp (i\varphi )=(z_2-z_1)/2
\end{equation}
and substitute $z=z_0+it\exp (i\varphi )$, then 
\begin{equation}
\label{eq53}J=\frac i{\exp [i\varphi (\alpha +\beta +1)]}\int_{-\infty
}^{+\infty }\frac{dt}{(r-it)^{\alpha +1}(r+it)^{\beta +1}}.
\end{equation}
The last integral can be found in tables (see e.g. \cite{tabi},p.301), using
(\ref{eq51}) we obtain 
\begin{equation}
\label{eq54}J=\frac{2i\pi }{(z_2-z_1)^{\alpha +\beta +1}}\frac{\Gamma
(\alpha +\beta +1)}{\Gamma (\alpha +1)\Gamma (\beta +1)}.
\end{equation}
Let us note that the opposite orientation of the line $K$ (i.e. point $z_2$
on left side in respect to the direction of $K,$ as the line $L$ in Fig.\ref
{fg6}$b$) should give result (\ref{eq54}) with opposite sign. Obviously the
integrals $I_3,I_4$ vanish and for $I_1,I_2$ we get 
\begin{equation}
\label{eq55}I_1=\frac{-2i\pi }{(x-z+2i\tau )^{\alpha +\beta +1}}\frac{\Gamma
(\alpha +\beta +1)}{\Gamma (\alpha +1)\Gamma (\beta +1)}
\end{equation}
\begin{equation}
\label{eq56}I_2=\frac{+2i\pi }{(x-z-2i\tau )^{\alpha +\beta +1}}\frac{\Gamma
(\alpha +\beta +1)}{\Gamma (\alpha +1)\Gamma (\beta +1)}
\end{equation}
and inserting into (\ref{eq49}) gives 
\begin{equation}
\label{eq57}I=\frac{(-1)^{\alpha +\beta +1}\Gamma (\alpha +\beta +1)}{2i\pi }%
\left( \frac 1{(x-z+2i\tau )^{\alpha +\beta +1}}-\frac 1{(x-z-2i\tau
)^{\alpha +\beta +1}}\right) .
\end{equation}
For $\tau \rightarrow 0$ this expression corresponds to $D^{\alpha +\beta
}(x-z)$ in (\ref{eq18}).

This result is formally correct, nevertheless its drawback is that it does
not reflect the correspondence of cut orientations in initial expressions (%
\ref{eq47}) and the final (\ref{eq57}). More rigorous discussion about the
cuts we postpone to Appendix, here give only result: the following
composition relation holds for operators (\ref{eq18}) with equally oriented
cuts 
\begin{equation}
\label{eq58}D_{\pm }^{\alpha +\beta }(x-z)=\int_{-\infty }^{+\infty }D_{\pm
}^\alpha (x-y)D_{\pm }^\beta (y-z)dy\qquad \alpha ,\beta \neq
-1,-2,-3...;\quad \alpha +\beta >-1. 
\end{equation}
Let us note, that validity of the relation (\ref{eq58}) can be verified in
the initial representation (\ref{eq7}) as well. Further, considering $D_{\pm
}^\alpha $ as a generalized function, then the condition $\alpha +\beta >-1$
can be omitted and the composition relation has the form 
\begin{equation}
\label{eq100}\int_{-\infty }^{+\infty }D_{\pm }^{\alpha +\beta }(x-\xi
)f(\xi )d\xi =\int_{-\infty }^{+\infty }\int_{-\infty }^{+\infty }D_{\pm
}^\alpha (x-y)D_{\pm }^\beta (y-\xi )f(\xi )d\xi dy,\quad \alpha ,\beta \neq
-1,-2,.. 
\end{equation}
for all functions analytic on the real axis and fulfilling 
\begin{equation}
\label{eq100+}\lim _{z\rightarrow \pm \infty }\frac{f(z)}{z^{\alpha +\beta }}%
=0. 
\end{equation}
Eq. (\ref{eq100}) follows from (\ref{eq58}) and relation $\frac d{dw}D_{\pm
}^{\gamma -1}(w)=D_{\pm }^\gamma (w).$ Repetitional integration by parts on
both sides can reduce sum $\alpha +\beta $ by any natural number, so in this
way validity of (\ref{eq100}) is proved.

All our previous considerations concerned acting of the operator ${\bf D}%
^\alpha $ on real axis. In the next we shall try to enlarge them on whole
complex plane.

\section{Fractional derivative in complex plane}

\begin{figure}
\begin{center}
\epsfig{file=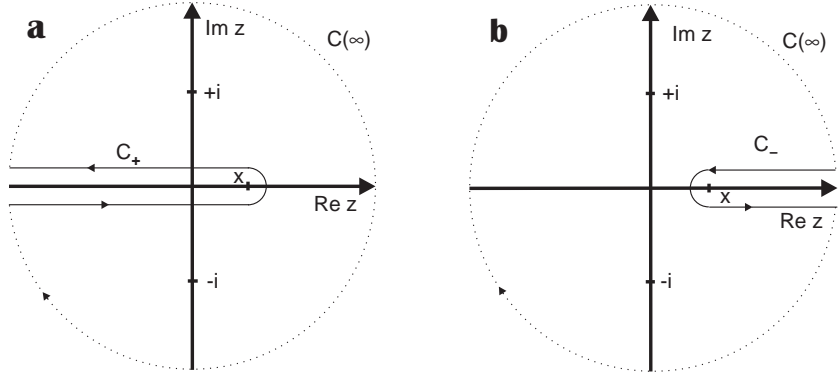,height=5cm}
\end{center}
\caption{Integration paths in Eq. (\ref{eq60}) closed in
infinity: $a)C_{+},\ b)C_{-}$.\label{fg8}} 
\end{figure}

\setcounter{equation}{0} First, let us illustrate the notion fractional
derivative introduced in previous part on a particular example. Take the
function 
\begin{equation}
\label{eq59a}f(x)=\frac 1{1+x^2}=-\frac 1{2i}\left( \frac 1{x+i}-\frac
1{x-i}\right) ,
\end{equation}
for which (\ref{eq33}) gives 
\begin{equation}
\label{eq60}f_{\pm }^\alpha (x)=\lim _{\epsilon \rightarrow 0+}\frac{%
(-1)^{\alpha +1}\Gamma (\alpha +1)}{2i\pi }\int_{C_{\pm }}\frac{dz}{%
(x-z)^{\alpha +1}(1+z^2)}.
\end{equation}
Obviously for $\alpha >-2$, the integration paths in Fig.\ref{fg5} closed by
arcs in infinity as shown in Fig.\ref{fg8}, are possible and 
\begin{equation}
\label{eq61}f_{\pm }^\alpha (x)=\frac{(-1)^{\alpha +1}\Gamma (\alpha +1)}{2i}%
\left( \frac 1{(x+i)^{\alpha +1}}-\frac 1{(x-i)^{\alpha +1}}\right) .
\end{equation}
Even though the result formally does not depend on given subscript (+ or $-$%
), it is necessary to take into account different cut orientation for $%
f_{+}^\alpha ,$ $f_{-}^\alpha .$ In accordance with the phase convention (%
\ref{eq6}), we take 
\begin{equation}
\label{eq62}
\begin{array}{lll}
\begin{array}{c}
(x+i)^{\alpha +1}= \\ 
(x-i)^{\alpha +1}=
\end{array}
& \left. 
\begin{array}{c}
R^{\alpha +1}\exp [i\varphi (\alpha +1)]\ \ \ \ \ \ \ \ \ \ \  \\ 
R^{\alpha +1}\exp [i(2\pi -\varphi )(\alpha +1)]\,\,
\end{array}
\right\}  & \qquad \qquad 
{\rm for\ }f_{+}^\alpha  \\ 
\begin{array}{c}
(x+i)^{\alpha +1}= \\ 
(x-i)^{\alpha +1}=
\end{array}
& \left. 
\begin{array}{c}
R^{\alpha +1}\exp [i\varphi (\alpha +1)]\ \ \ \ \ \ \ \ \ \ \  \\ 
R^{\alpha +1}\exp [-i\varphi (\alpha +1)]\ \ \ \ \ \ \ \ \ 
\end{array}
\right\}  & \qquad \qquad {\rm for\ }f_{-}^\alpha ,
\end{array}
\end{equation}
where 
\begin{equation}
\label{eq62+}R=\sqrt{1+x^2},\qquad \qquad \varphi =\frac \pi 2-\arcsin \frac
xR.
\end{equation}
After inserting into (\ref{eq61}) and a simple rearrangement, we obtain 
\begin{equation}
\label{eq63}f_{\pm }^\alpha (x)=\left( \frac{\pm 1}{\sqrt{1+x^2}}\right)
^{\alpha +1}\Gamma (\alpha +1)\sin \left[ \left( \arcsin \frac x{\sqrt{1+x^2}%
}\pm \frac \pi 2\right) \left( \alpha +1\right) \right] .
\end{equation}

\noindent {\it Remark: }From the last relation it can be easily shown that
e.g.

$a)$ for $\alpha $ integer, $\alpha =n>-1$, it holds $f_{+}^n(x)=f_{-}^n(x)%
\equiv f^n(x),$ where $f^n$ is ordinary {\it n-fold} derivative

$b)$ $f_{-}^\alpha (x)=(-1)^\alpha f_{+}^\alpha (-x)$

$c)$ for $\alpha \rightarrow -1$ we get a primitive function to $f:\quad
f_{\pm }^\alpha (x)\rightarrow \arctan (x)\pm \pi /2.$ The same results can
be obtained also from (\ref{eq46}) for $n=-1,\quad \Delta \alpha =0.$ \\

Using formula (\ref{eq33}) we shall try to generalize the operator of
fractional derivative on real axis to the whole complex plane. For this
operator we shall demand again: 
\begin{equation}
\label{eq64a}{\bf D}^\alpha f(z)=\frac{d^nf}{dz^n}\qquad {\rm for\ }\alpha
=n\geq 0,
\end{equation}
\begin{equation}
\label{eq64b}{\bf D}^\alpha \,\circ {\bf D}^\beta ={\bf D}^{\alpha +\beta }.
\end{equation}
For making the generalization more transparent, we shall do it in several
steps.

\begin{figure}
\begin{center}
\epsfig{file=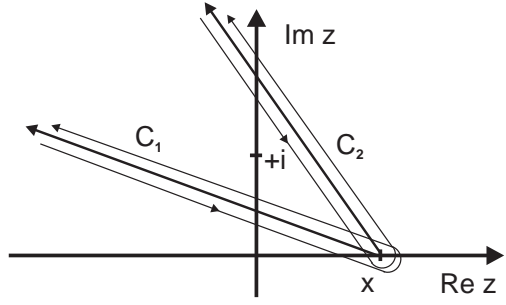,height=4cm}
\end{center}
\caption{Tentative modification of integration paths in Eq. (\ref{eq60}) 
and Fig.\ref{fg8}.\label{fg9}} 
\end{figure}

\begin{figure}
\begin{center}
\epsfig{file=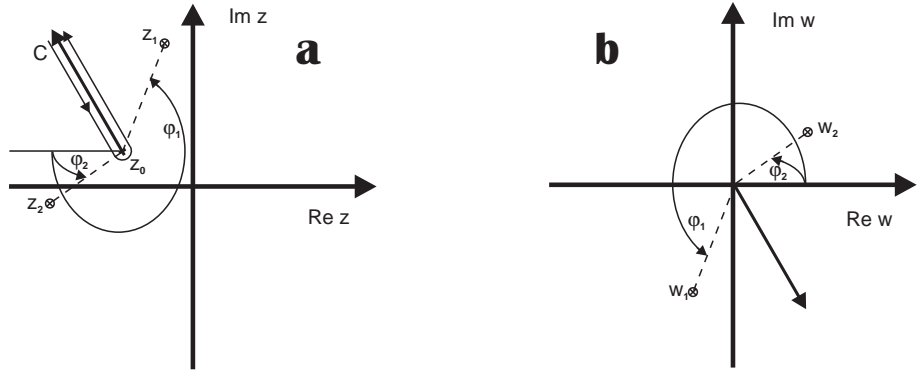,height=5cm}
\end{center}
\caption{$a)$Integration paths in Eq. (\ref{eq33}) for the function 
(\ref{eq65a}). $\varphi _{1},\ \varphi _{2}$ are phases corresponding to terms 
$w_k^{\alpha +1}=(z_0-z_k)^{\alpha +1}$
in the result of integration (\ref{eq66}). $b)$ The same phases represented
for variable $w$.\label{fg10}} 
\end{figure}

1) Let us go back to Fig.\ref{fg8} and put a question what will change in
results (\ref{eq61})$-$(\ref{eq63}), if the cut (and integration path
correspondingly) would be oriented otherwise than along real axis, but e.g.
as shown in Fig.\ref{fg9}. Obviously for Eq. (\ref{eq61}) nothing will 
change, but the form of Eqs. (\ref{eq62}),(\ref{eq63}) will depend on the
mutual position of the cut and both poles. Because we have accepted phase
convention (\ref{eq6}) only for cuts on real axis, it is now necessary to
make a more general consideration to determine the phases of both terms in (%
\ref{eq61}). Instead of (\ref{eq59a}) let us take a complex function 
\begin{equation}
\label{eq65a}f(z)=\frac{a_1}{z-z_1}+\frac{a_2}{z-z_2} 
\end{equation}
and integration paths $C$ displayed in Fig.\ref{fg10}$a.$ The corresponding 
integral will be 
\begin{equation}
\label{eq65b}f^\alpha (z_0)=-(-1)^{\alpha +1}\Gamma (\alpha +1)\left( \frac{%
a_1}{(z_0-z_1)^{\alpha +1}}+\frac{a_2}{(z_0-z_2)^{\alpha +1}}\right) . 
\end{equation}
The path $C$ passes about the cut of function $1/w^{\alpha
+1}=1/(z_0-z)^{\alpha +1}$ in variable $z$, i.e. cut of the function $%
1/w^{\alpha +1}$ is oriented in opposite direction, see Fig.\ref{fg10}$b.$
The phases of $w_1^{\alpha +1}=(z_0-z_1)^{\alpha +1}$ and $w_2^{\alpha
+1}=(z_0-z_2)^{\alpha +1}$ must be fixed in respect to this cut. Fig.\ref
{fg10}$b$ prompts the following rule.

\noindent {\bf Rule 1: }{\it Phase of complex variable $w$ is given by the
angle $\varphi $ of arc leading from positive real half axis and measured in
positive direction (against clockwise sense) to the point $w$ and if the arc
intersects the cut, $\varphi $ is reduced by 2$\pi .$}

\noindent
This rule can be applied also directly for situation in Fig.\ref{fg10}$a$
for fixing phase of $(z_0-z)^{\alpha +1}$. The only modification is that $%
\varphi $ is measured from half line $(z_0,z_0-\infty ).$ Therefore angles $%
\varphi _1,\varphi _2$ in Fig.\ref{fg10} fix phases in (\ref{eq65b}) 
\begin{equation}
\label{eq66}f^\alpha (z_0)=(-1)^\alpha \Gamma (\alpha +1)\left( \frac{%
a_1\exp (-i\varphi _1[\alpha +1])}{\mid (z_0-z_1)^{\alpha +1}\mid }+\frac{%
a_2\exp (-i\varphi _2[\alpha +1])}{\mid (z_0-z_2)^{\alpha +1}\mid }\right) . 
\end{equation}
Now the introduced rule can be applied also for integration on paths $%
C_1(C_2)$ in Fig.\ref{fg9}. Doing this, we shall get the same result as in
the case of integration on paths $C_{+}(C_{-})$ in Fig.\ref{fg8}. Therefore
depending on the position of chosen derivative cut in respect to the both
poles, the function ( \ref{eq59a}) has (up to factor (\ref{eq8})) in a given
point $x$ two different values of fractional derivative given by Eq. (\ref
{eq63}).

\begin{figure}
\begin{center}
\epsfig{file=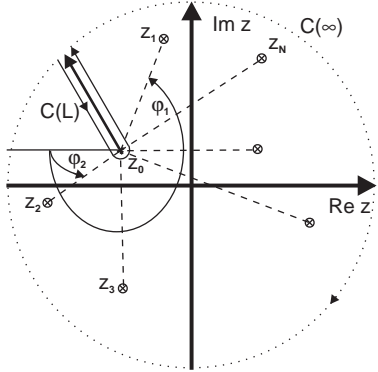,height=5cm}
\end{center}
\caption{Integration paths in Eq. (\ref{eq33}) with the function (\ref{eq67})
and phases $\varphi _k$ appearing in the result of integration (\ref{eq68}).
\label{fg11}}
\end{figure}

2)We shall now generalize the prescription for calculating of fractional
derivative of function having the form (\ref{eq65a}) for functions having a
finite number of poles 
\begin{equation}
\label{eq67}h(z)=\sum_{k=1}^N\frac{a_k}{(z-z_k)^{n_k+1}},\qquad \qquad
n_k\geq 0.
\end{equation}
Let the derivative cut be given by half line $L\equiv [z_0,\ z_0+\exp
(i\theta )\infty ]$ which does not go through any of poles $z_k.$ Then 
$$
h^\alpha (z_0)=\lim _{\epsilon \rightarrow 0+}\frac{(-1)^{\alpha +1}\Gamma
(\alpha +1)}{2i\pi }\sum_{k=1}^N\int_{C(L)}\frac{a_kdz}{(z_0-z)^{\alpha
+1}(z-z_k)^{n_k+1}}\quad \ \  
$$
\begin{equation}
\label{eq68}=(-1)^\alpha \sum_{k=1}^N\frac{\Gamma (\alpha +n_k+1)}{\Gamma
(n_k+1)}\frac{a_k\exp (-i\varphi _k[\alpha +1])}{\mid z_0-z_k\mid ^{\alpha
+1}(z_0-z_k)^{n_k}},
\end{equation}
integration path $C(L)$ is shown in Fig.\ref{fg11}. 
Angles $\varphi _k$ are calculated using the {\bf Rule 1}, therefore it is
obvious, that the function $h^\alpha (z_0)$ in (\ref{eq68}) can have,
according to the cut orientation, as much values, as much different poles $%
z_k$ the function (\ref{eq67}) has.

\begin{figure}
\begin{center}
\epsfig{file=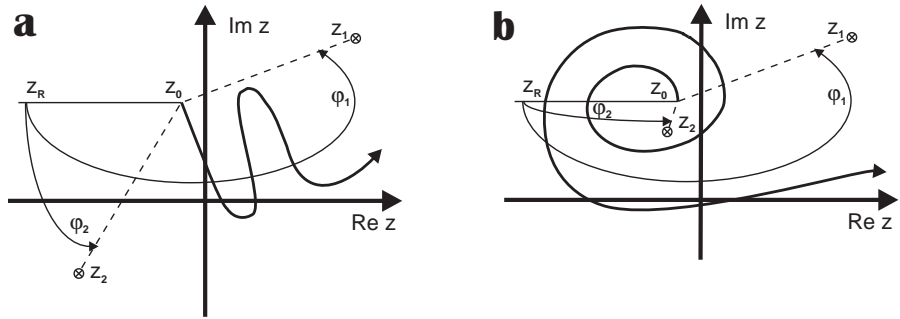,width=12cm}
\end{center}
\caption{$a),\ b)$ Examples of the curvilinear cuts generating integration
paths for Eq. (\ref{eq33}) with the function (\ref{eq67}). Other symbols are
defined in {\bf Rule 2}.\label{fg12}}
\end{figure}

3) Now let us consider functions (\ref{eq67}),(\ref{eq68}) in the case, when
the cut is not a line, but some general curve connecting the points $z_0$
and $\exp (i\theta )\infty .$ Then the result (\ref{eq68}) will be formally
the same, only the prescription for fixing of angles $\varphi _k$ must be
modified. General cases are illustrated in Fig.\ref{fg12}{\it a,b}, where it 
is apparent that the arc, which 'measures' angle, may intersect the cut
several times. The consistent assigning of angles $\varphi _k$ corresponding
to points $z_k$ in Eq. (\ref{eq68}) can be ensured by the following
prescription.

\noindent {\bf Rule 2: }{\it Let us choose on the half line $(z_0,z_0-\infty
)$ any reference point $z_R.$ Angle $\varphi _k$ is given by angle $%
z_Rz_0z_k $ measured in positive sense and then for each intersection with
the cut is corrected as follows. Superpose palm of right or left hand at an
intersection in such a way, that fingers lead from $z_R$ towards z$_k$ and
thumb leads in direction of the cut from branching point $z_0$. If this
condition is met by right (left) hand, $\varphi _k$ will be enhanced
(reduced) by 2$\pi $.}

{\it Remark: }It is substantial for all $z_k$ to choose one common reference
point $z_R.$ A shift of this point results at most in only equal shifts of
all angles $\varphi _k\rightarrow \varphi _k+2n\pi $, which do not change
Eq. (\ref{eq68}). \\

Therefore for any curvilinear cut the relation (\ref{eq68}) can be written 
\begin{equation}
\label{eq70}h^\alpha (z_0)=(-1)^\alpha \sum_{k=1}^N\frac{\Gamma (\alpha
+n_k+1)}{\Gamma (n_k+1)}\frac{a_k\exp (-i[\varphi _k+2m_k\pi ][\alpha +1])}{%
\mid z_0-z_k\mid ^{\alpha +1}(z_0-z_k)^{n_k}}.
\end{equation}
The set of integer numbers (or more exactly their differences) characterizes
the way how the curvilinear cut passes among the poles $z_k$. If we accept
curvilinear cuts for derivative of the function (\ref{eq59a}), then using (%
\ref{eq70}) one can obtain 
\begin{equation}
\label{eq72}f_{}^\alpha (x)=\frac{(-1)^{k(\alpha +1)}\Gamma (\alpha +1)}{%
\left( \sqrt{1+x^2}\right) ^{\alpha +1}}\sin \left[ \left( \arcsin \frac x{
\sqrt{1+x^2}}\,+\,(2k+1)\frac \pi 2\right) \left( \alpha +1\right) \right] ,
\end{equation}
where $k=m_2-m_1$. Let us note that 
\begin{equation}
\label{eq73}\lim _{\alpha \rightarrow -1}f^\alpha (x)=\arctan
(x)+(2k+1)\frac \pi 2,
\end{equation}
i.e. we obtain the infinite (but countable) set of primitive functions for
the function (\ref{eq59a}).

\begin{figure}
\begin{center}
\epsfig{file=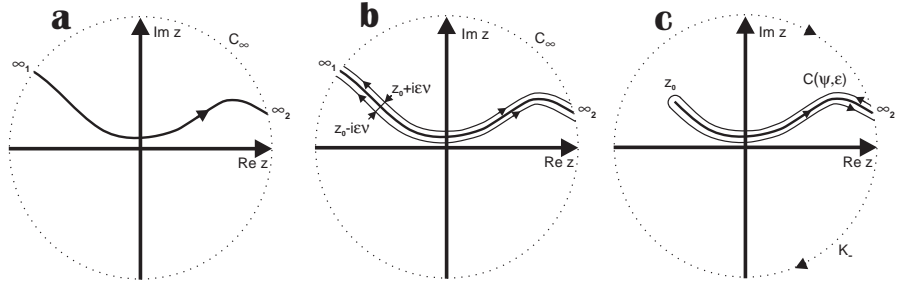,width=12cm}
\end{center}
\caption{$a)$Integration paths for operator (\ref{eq76}) in complex plane. 
$b)$ Corresponding curvilinear cuts in Eq. (\ref{eq76}).
$c)$ Corresponding integration path in Eq. (\ref{eq79}).\label{fg15}}
\end{figure}

4) So far we have considered only analytic functions of the form (\ref{eq67}%
), for which the corresponding integral on whole circle $C(\infty )$
vanishes. Now we are going to consider the general case, when this integral
vanishes only on a part of the circle. But first, let us go back to the
operator (\ref{eq18}) acting on analytic functions on real axis (or the part
of this axis) and try to generalize it for analytic functions on a curve in
complex plane connecting a pair of points on the circle $C(\infty ),$ see
Fig.\ref{fg15}$a$. Let this curve be given as a continuous complex function
of 
real parameter $t\in (-\infty ,+\infty )$ 
\begin{equation}
\label{eq74}z=x(t)+iy(t)\equiv \psi (t),\qquad \psi (-\infty )=\infty
_1,\qquad \psi (+\infty )=\infty _2,
\end{equation}
then its derivative 
\begin{equation}
\label{eq75}\psi ^{\prime }(t)=\frac{dx}{dt}+i\frac{dy}{dt}
\end{equation}
determines in complex plane the vector tangent to curve $\psi $ and oriented
in direction of increasing $t.$ The function 
\begin{equation}
\label{eq75+}\nu (t)=\frac{\psi ^{\prime }(t)}{\sqrt{\psi ^{\prime }(t)\cdot
\psi ^{*\prime }(t)}}=\exp [i\omega (t)]
\end{equation}
represents this vector in its normalized value ($\omega (t)$ is phase of
this vector) and the function $i\nu (t)$ normalized vector perpendicular to $%
\psi $ and oriented left in respect to the course of $\psi $.

Now let us define the function of the complex variable $z\in \psi $ with
complex parameters $z_0\in \psi $ and $\alpha \neq -1,-2,-3,...$ 
\begin{equation}
\label{eq76}D_\psi ^\alpha (z_0-z)=\lim _{\epsilon \rightarrow 0+}\frac{%
(-1)^{\alpha +1}\Gamma (\alpha +1)}{2i\pi }\left( \frac 1{(z_0-z+i\epsilon
\nu )^{\alpha +1}}-\frac 1{(z_0-z-i\epsilon \nu )^{\alpha +1}}\right) 
\end{equation}
having the curvilinear cuts (for $\alpha \neq 0,1,2,..)$ corresponding to
the both terms coming out the points $z_0\pm i\epsilon \nu $ and going
jointly along the curve $\psi $ to points $\infty _1$ or $\infty _2$, see
Fig.\ref{fg15}$b.$ So, in contradistinction to linear cuts the form of the
cut of $w^{\alpha +1}=(z_0-z)^{\alpha +1}$ does depend on position of point $%
z_0$ on $\psi .$

Next, using this function we define the operator 
\begin{equation}
\label{eq77}{\bf D}_\psi ^\alpha \ f=f_{\psi \pm }^\alpha (z_0)=\int_\psi
D_\psi ^\alpha (z_0-z)f(z)dz=\int_{-\infty }^{+\infty }D_\psi ^\alpha (\psi
(t_0)-\psi (t))f(\psi (t))\psi ^{\prime }(t)dt 
\end{equation}
acting on the functions analytic on the curve $\psi $ for which the integral
converges.

\noindent {\it Remark:}{\bf \ }Let us note the integral (\ref{eq77}) depends
on choice of cut end-point ($\infty _1$ or $\infty _2)$ but does not depend
on in which direction ($\infty _1$ $\rightarrow $ $\infty _2\ $ or\ $\infty
_1$ $\leftarrow $ $\infty _2)$ integration is done. That is a consequence of
the fact that change of integration direction $\psi (t)\rightarrow \psi
(-t),\ dz\rightarrow -dz$ implies the change $\nu (t)\rightarrow -\nu (t),$
which implies the change $D_\psi ^\alpha (z_0-z)\rightarrow -D_\psi ^\alpha
(z_0-z),$ therefore the integral does not change.\\

If we calculate (\ref{eq77}) as the difference 
\begin{equation}
\label{eq78}{\bf D}_\psi ^\alpha \ f=\lim _{\epsilon \rightarrow 0+}\frac{%
(-1)^{\alpha +1}\Gamma (\alpha +1)}{2i\pi }\left( \int_\psi \frac{f(z)dz}{%
(z_0-z+i\epsilon \nu )^{\alpha +1}}-\int_\psi \frac{f(z)dz}{(z_0-z-i\epsilon
\nu )^{\alpha +1}}\right) , 
\end{equation}
then this difference can be expressed as 
\begin{equation}
\label{eq79}\int_\psi D_\psi ^\alpha (z_0-z)f(z)dz=\lim _{\epsilon
\rightarrow 0+}\frac{(-1)^{\alpha +1}\Gamma (\alpha +1)}{2i\pi }\int_{C(\psi
,\epsilon )}\frac{f(z)dz}{(z_0-z)^{\alpha +1}}, 
\end{equation}
where the path in 'distance' $\epsilon $ passes about the cut coming out
branching point $z_0$ on the curve $\psi $ to infinity, see Fig.\ref{fg15}$%
c. $

To label the integration cuts ending either at $\infty _1$ or $\infty _2$ we
can accept the following convention. Let $\infty _1=\exp (i\theta _1)\infty $
and $\infty _2=\exp (i\theta _2)\infty ,$ then for 
\begin{equation}
\label{eq80}
\begin{array}{llll}
\cos \theta _1\neq \cos \theta _2 & \left\{ 
\begin{array}{c}
\cos \theta _1<\cos \theta _2 \\ 
\cos \theta _1>\cos \theta _2
\end{array}
\right\}  & {\rm we{\rm \ }define} & \left\{ 
\begin{array}{c}
\psi _{+}=(z_0,\infty _1),\quad \psi _{-}=(z_0,\infty _2) \\ 
\psi _{+}=(z_0,\infty _2),\quad \psi _{-}=(z_0,\infty _1)
\end{array}
\right\}  \\ 
\cos \theta _1=\cos \theta _2 & \left\{ 
\begin{array}{c}
\ \sin \theta _1<\sin \theta _2\  \\ 
\sin \theta _1>\sin \theta _2
\end{array}
\right\}  & {\rm we{\rm \ }define} & \left\{ 
\begin{array}{c}
\psi _{+}=(z_0,\infty _1),\quad \psi _{-}=(z_0,\infty _2) \\ 
\psi _{+}=(z_0,\infty _2),\quad \psi _{-}=(z_0,\infty _1)
\end{array}
\right\} 
\end{array}
\end{equation}
and correspondingly we index related symbols ${\bf D}_{\psi \pm }^\alpha ,\
D_{\psi \pm }^\alpha ,\ \ f_{\psi \pm }^\alpha ,\ C_{\pm }(\psi ,\epsilon ).$

Let us note the definition of the operator ${\bf D}_\psi ^\alpha $ by Eqs. (%
\ref{eq76}), (\ref{eq77}) ensure that the both curves $C_{\pm }(\psi
,\epsilon )$ are oriented in such a way that after their closing by the
circle $C(\infty )$ see Fig.\ref{fg15}$c$, there arises closed curve having 
{\it always}{\bf \ }clockwise orientation. We denote these curves as $K_{\pm
}(\psi ,\epsilon ).$

Now, after all these preparing steps we can come to the formulation and
proof of the following theorem.

{\bf Theorem}

\noindent
{\it Assumptions:}

\begin{figure}
\begin{center}
\epsfig{file=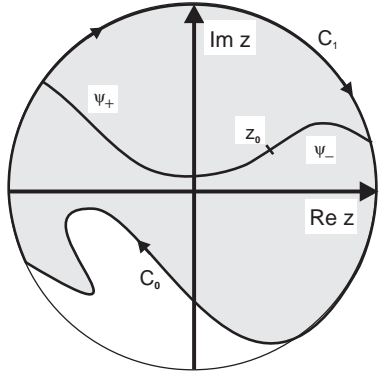,height=5cm}
\end{center}
\caption{Domain $G$ (grey area) and the curves defined in assumptions 
of the 
{\bf Theorem}.\label{fg16}} 
\end{figure}

$i)$ $G$ is the domain in complex plane containing the part (or parts) $C_1$ 
$\equiv G\cap C(\infty )$ of the circle $C(\infty ).$ The curve $C_1$ forms
one part, the curve $C_0$ is the second part and closed curve $C_G\equiv
C_0\cup C_1$ constitutes the complete domain boundary (Fig. \ref{fg16}). 

$ii)$ $\psi $ is a curve defined in (\ref{eq74}) and $\psi \subset G,\ \psi
\cap C_0=\emptyset .$

$iii)$ Function $f(z)$ is analytic in $G$ except for the poles $\ z_k\notin
\psi ,\ z_k\notin C_G,\ k=1...N$ and for a given $\alpha \equiv \alpha
_1+i\alpha _2\neq -1,-2,-3...\ $ and any $z_1\in C_1\ $meets the condition 
\begin{equation}
\label{eq81}\lim _{z\rightarrow z_1}\frac{f(z)}{z^{\alpha _1}}=0. 
\end{equation}

\noindent
{\it Statements:}

1. Operation 
\begin{equation}
\label{eq83}{\bf D}_{\psi \pm }^\alpha f=f_{\psi \pm }^\alpha
(z_0)=\int_\psi D_{\psi \pm }^\alpha (z_0-z)f(z)dz,\qquad z_0\in \psi 
\end{equation}
is a fractional derivative satisfying the conditions (\ref{eq64a}), (\ref
{eq64b}). In particular that means the operation does not depend on $\psi $
if $\alpha =n\geq 0$ is integer.

2. $f_{\psi \pm }^\alpha $ is given also equivalently by 
\begin{equation}
\label{eq84}f_{\psi \pm }^\alpha (z_0)=\lim _{\epsilon \rightarrow 0+}\frac{%
(-1)^{\alpha +1}\Gamma (\alpha +1)}{2i\pi }\int_{C_{\pm }(\psi ,\epsilon )} 
\frac{f(z)dz}{(z_0-z)^{\alpha +1}}. 
\end{equation}

3. For $N=0$ the value $f_{\psi \pm }^\alpha (z_0)$ does not depend on $\psi
,$ i.e. up to factor $(-1)^\alpha $ is uniquely defined.

4. For $N\geq 1$ and $\alpha $ non integer there remain cases:

4.1. $\alpha _1$ is rational and $\alpha _2=0$, then $f_{\psi \pm }^\alpha
(z_0)$ does depend on $\psi $, but the set of values is finite.

4.2. $\alpha _1$ is irrational or $\alpha _2\neq 0$, then $f_{\psi \pm
}^\alpha (z_0)$ does depend on $\psi $ as well and in general the set of
values is infinite, but countable.\\

\noindent
{\it Proof:}

\noindent
We shall start from the statement 2. Obviously, its validity follows from
Eqs. (\ref{eq78}), (\ref{eq79}).

Now let us consider the statement 4. If we express the function $f(z)$ as
the sum $f(z)=g(z)+h(z)$, where $g(z)$ is analytic in $G$ and $h(z)$ has the
form (\ref{eq67}), then obviously 
$$
f_{\psi \pm }^\alpha (z_0)=\lim _{\epsilon \rightarrow 0+}\frac{(-1)^{\alpha
+1}\Gamma (\alpha +1)}{2i\pi }\int_{C_{\pm }(\psi ,\epsilon )}\frac{f(z)dz}{%
(z_0-z)^{\alpha +1}} 
$$
\begin{equation}
\label{eq93}=\lim _{\epsilon \rightarrow 0+}\frac{(-1)^{\alpha +1}\Gamma
(\alpha +1)}{2i\pi }\int_{K_{\pm }(\psi ,\epsilon )}\frac{f(z)dz}{%
(z_0-z)^{\alpha +1}}-\frac{(-1)^{\alpha +1}\Gamma (\alpha +1)}{2i\pi }%
\int_{C_0}\frac{f(z)dz}{(z_0-z)^{\alpha +1}} 
\end{equation}
$$
=(-1)^\alpha \left( \sum_{k=1}^N\frac{\Gamma (\alpha +n_k+1)}{\Gamma (n_k+1)}
\frac{a_k\exp (-i[\varphi _k+2m_k\pi ][\alpha +1])}{\mid z_0-z_k\mid
^{\alpha +1}(z_0-z_k)^{n_k}}+I(z_0)\right) , 
$$
where we have denoted 
\begin{equation}
\label{eq82}I(z_0)=\frac{\Gamma (\alpha +1)}{2i\pi }\int_{C_0}\frac{f(z)dz}{%
(z_0-z)^{\alpha +1}}. 
\end{equation}
The angles $\varphi _k+2m_k\pi $ corresponding to the poles $z_k$ are
evaluated according to {\bf Rule 2. }In the function $I(z_0)$ the
orientation of the path $C_0$ is accordant with $K_{\pm }$ and phase of $%
(z_0-z)^{\alpha +1}$ must be also according to the {\bf Rule 2} related to
the same reference point $z_R$. Now, if $\alpha _1=p/q$ ($p,q$ are not
commensurable) and $\alpha _2=0$, then each term in the last sum has the
same value also for $m_k^{\prime }=m_k+q$, i.e. depending on the form of the
cut and the corresponding set $\{m_k,1\leq k\leq N\}$ the expression (\ref
{eq93}) has only {\it finite} number of values for some $z_0$. On the other
hand for $\alpha $ irrational different sets $\{m_k\}$ give different values
of the sum, therefore in general the number of values is {\it infinite.} For 
$\alpha _2\neq 0$ multi-value factors in the sum (\ref{eq93}) are expanded 
\begin{equation}
\label{eq99}1^\alpha =\exp (2i\pi m\alpha _1)\cdot \exp (-2\pi m\alpha _2) 
\end{equation}
i.e. $m\alpha _1,m\alpha _2$ determine phase and scale of individual terms
in the sum. Obviously for a complex $\alpha $ the set of values $f^\alpha
(z_0)$ depending on $\psi $ is in general {\it infinite}, similarly as for $%
\alpha $ irrational. Therefore the statements 4.1, 4.2 are proved.

For $N=0$ Eq. (\ref{eq93}) is simplified: 
\begin{equation}
\label{eq94}f_{\psi \pm }^\alpha (z_0)=(-1)^\alpha I(z_0) 
\end{equation}
and trueness of statement 3. is evident.

Finally, let us consider the statement 1. Validity of the condition (\ref
{eq64a}) follows from the statement 2. For $\alpha $ integer the path $%
C_{\pm }(\psi ,\epsilon )$ can be closed around the pole $z_0$ (having {\it %
positive} orientation) and we get the Cauchy integral 
\begin{equation}
\label{eq95}\frac{n!}{2i\pi }\int_{C(z_0)}\frac{f(z)dz}{(z-z_0)^{n+1}}%
=f^n(z_0). 
\end{equation}
The condition (\ref{eq64b}) requires validity 
\begin{equation}
\label{eq96}\int_\psi D_{\psi \pm }^{\alpha +\beta }(x-\xi )f(\xi )d\xi
=\int_\psi \int_\psi D_{\psi \pm }^\alpha (x-y)D_{\psi \pm }^\beta (y-\xi
)f(\xi )d\xi dy. 
\end{equation}
This relation follows from (\ref{eq77}) and (\ref{eq100}) in which the
substitutions $\xi \rightarrow \psi (\xi )$ and $y\rightarrow \psi (y)$ are
applied. So the whole proof is completed.

\section{Discussion}

\setcounter{equation}{0}

\subsection{Some remarks regarding the theorem}

Now, let us look into the statements of the {\bf Theorem }more
comprehensively. The function $I(z_0)$ in (\ref{eq93}) can be expressed 
$$
I(z_0)=\frac{\Gamma (\alpha +1)}{2i\pi }\int_{C_0}\frac{f(z)dz}{%
(z_0-z)^{\alpha +1}}= 
$$
\begin{equation}
\label{eq97}=\frac{\Gamma (\alpha +1)\exp (-i2m\pi [\alpha +1])}{2i\pi }%
\int_{C_0}\frac{f(z)\exp [-i\varphi (z)(\alpha +1)]dz}{\left|
(z_0-z)^{\alpha +1}\right| }. 
\end{equation}
Since the cut $\psi $ does not intersect the curve $C_0$ the phases of all
points on $C_0$ are 'corrected' by the same factor standing ahead of the
integral. Now it is obvious that for $\alpha =p/q$ the function (\ref{eq93})
can have at most $q^{N+1}$ different values, including multi-value factor $%
(-1)^\alpha $.

Actually, the function (\ref{eq93}) can be written as 
\begin{equation}
\label{eq98}f^\alpha (z)=(-1)^\alpha \left( \sum_{k=1}^N\frac{\Gamma (\alpha
+n_k+1)}{\Gamma (n_k+1)}\frac{a_k}{(z-z_k)^{n_k+\alpha +1}}+\frac{\Gamma
(\alpha +1)}{2i\pi }\int_{C_0}\frac{f(\xi )d\xi }{(z-\xi )^{\alpha +1}}%
\right) 
\end{equation}
and considered a multi-value function with the value at point $z$ depending
on the choice of cut $\psi $ in (\ref{eq83}). At given point $z$ there are
equivalent any two cuts for which the region closed by them and the curve $%
C_1$ does not contain any pole $z_k$. Obviously, except for the points $z_k$
the function $f^\alpha (z)$ is (as original function $f$) analytic in the
domain $G.$ Moreover, in the region of $\alpha $ in which $f^\alpha ={\bf D}%
^\alpha f$ exists, this function is apparently analytic also in respect to $%
\alpha .$

Special case takes place when $\alpha $ is a negative integer. Then due to
the singularity $\Gamma (-n)$ the kernel (\ref{eq76}) loses the sense.
Nevertheless, assuming in (\ref{eq76}) initially $\alpha \neq -n$ one can
proceed to representation (\ref{eq84}), then making the limit $\epsilon
\rightarrow 0$ (like Eq. (\ref{eq38})) gives 
\begin{equation}
\label{eq101}f_{\psi \pm }^\alpha (z_0)=-\lim _{\epsilon \rightarrow 0+} 
\frac{\Gamma (\alpha +1)\sin ([\alpha +1]\pi )}\pi \left( \int_{\psi \pm } 
\frac{f(z+\nu _0\epsilon )dz}{(z_0-z-\nu _0\epsilon )^{\alpha +1}}+\frac{%
f(z_0)}{\alpha (-\nu _0\epsilon )^\alpha }\right) , 
\end{equation}
where $\nu _0$ is given by Eq. (\ref{eq75+}) and represents direction of the
cut at $z_0$ (orientation is assumed from $z_0$ to infinity). This formula
already makes sense for $\alpha \ =-n$ (if condition (\ref{eq81}) holds).
Using relation (\ref{A7}) gives 
\begin{equation}
\label{eq102}f_{\psi \pm }^{-n}(z_0)=-\frac 1{(n-1)!}\int_{\psi \pm
}(z_0-z)^{n-1}f(z)dz, 
\end{equation}
which is a modification of the known formula (\ref{i1}) for {\it n-fold}
primitive function. The set of values $f^{-n}(z_0)$ depends on $\psi $ as
follows. If $f(z)$ has no poles in the domain $G$, then any two integration
paths $\psi _1,\psi _2$ in (\ref{eq102}) can be connected by a fragment of $%
C_1$, in this way there arises closed curve $C$ and it holds 
\begin{equation}
\label{eq103}0=-\frac 1{(n-1)!}\int_C(z_0-z)^{n-1}f(z)dz=f_{\psi
_2}^{-n}(z_0)-f_{\psi _1}^{-n}(z_0) 
\end{equation}
i.e. $f^{-n}(z_0)$ is determined uniquely. Now, let us suppose $f(z)$ has
inside curve $C$ just one pole, for $z\rightarrow z_p$ 
\begin{equation}
\label{eq104}f(z)\rightarrow \frac{a_p}{(z-z_p)^{n_p+1}}. 
\end{equation}
Then obviously any curve $C$ in the integral (\ref{eq103}) can be always
substituted by a couple of curves $K_0$, each of them is closed in the phase
range $\left\langle 0,2\pi \right\langle $ and having the pole $z_p$ inside
(e.g. circles centered at $z_p$). Instead of (\ref{eq103}) one get the
difference 
\begin{equation}
\label{eq105}\Delta _p\equiv f_{\psi _2}^{-n}(z_0)-f_{\psi _1}^{-n}(z_0)= 
\frac{ma_p}{(n-1)!}\int_{K_0}\frac{(z_0-z)^{n-1}}{(z-z_p)^{n_p+1}}dz, 
\end{equation}
where $m$ is the integer depending on the shape of the curve $C$ ($m$
represents the number of 'twists' on the $C$). Using the Cauchy formula, the
last equation gives 
\begin{equation}
\label{eq106}\Delta _p=\left\{ 
\begin{array}{cc}
0 & n_p\geq n \\ 
\frac{
\begin{array}{c}
2i\pi ma_p 
\end{array}
}{
\begin{array}{c}
(n-n_p-1)!n_p! 
\end{array}
}(z_0-z_p)^{n-n_p-1} & n_p<n 
\end{array}
\right. 
\end{equation}
For more poles the all corresponding terms (\ref{eq106}) are simply added.
For example integration constants in (\ref{eq73}) representing the case $%
n=1,\ n_p=0$ fulfill (\ref{eq106}). The composition relation for $\alpha
=-n,\ \beta <0$ and $\psi \equiv ${\it real axis} is proved in Appendix (Eq.
(\ref{A17})) and apparently can be transformed to any curvilinear path $\psi
. $

Let us note, for {\it $\alpha $ negative integer} only the representation of 
${\bf D}^\alpha $ given by Eq. (\ref{eq101}) (or (\ref{eq102})) makes sense
and conversely for {\it $\alpha $ non negative integer} only the
representation given by Eq. (\ref{eq83}) with kernel (\ref{eq76}) (or
equivalently by Eq. (\ref{eq84})) is well defined. For any {\it $\alpha $
complex but non integer} both representations are well defined. These two
representations differ in the corresponding integration paths:

$i)$ Eq. (\ref{eq84}) - integration on some curve enveloping the cut, the
path can be closed.

$ii)$ Eq. (\ref{eq101}) - integration on the cut itself, the path cannot be
closed.

\noindent
For $\alpha $ {\it integer} it is specific, that the cuts disappear.

Perhaps most restrictive assumption in the {\bf Theorem }is the condition (%
\ref{eq81}). The question, if one can in some consistent way apply ${\bf D}%
^\alpha $ to the functions not obeying this condition in any part of $%
C(\infty )$ (and simultaneously having the ordinary derivatives - primitive
functions) requires a further study. Obviously one possible way is to
consider such functions as generalized functions as well.

\subsection{Concluding remarks}

Have we said something new? First let us show what is not new. Apparently:

{\it a) }The content of Eq. (\ref{eq46}) is almost identical with Liouville
definition of right and left -handed fractional derivatives in \cite{samko},
p.95. The only distinction is in phase of $f_{-}^\alpha $ since in the cited
definition for the real functions the real value of the derivatives is ad
hoc postulated.

{\it b) }Eq. (\ref{eq84}) is the Cauchy type integral which ordinarily
serves as one of possible starting points for the definition of the
fractional derivative in complex plane, see \cite{samko},p.415. In our
approach the integration path is uniquely defined by chosen cut. \\

\noindent
On the other hand, the new seems be the following:

1) The general form of the kernel (\ref{eq76}) from which the both above
mentioned formulae follow.

2) The construction based on the integration paths enveloping curvilinear
cuts, which in the result allows to identify fractional derivative -
integral with the multi-valued function and to determine how the number of
values depends on the derivative order type and the number of poles which
the given function has in the considered region.\\
\newpage
\appendix


{\bf Appendix: The correspondence of cuts in the composition relation}

\setcounter{equation}{0} 
\renewcommand{\theequation}{A.\arabic{equation}} Similarly as in (%
\ref{eq48}) we define 
\begin{equation}
\label{A1} 
\begin{array}{ll}
\begin{array}{c}
h^{\pm }(\gamma ,w)= 
\frac{
\begin{array}{c}
1 
\end{array}
}{
\begin{array}{c}
(w+i\tau )^\gamma 
\end{array}
} \\ h_{\pm }(\gamma ,w)=-\frac{
\begin{array}{c}
1 
\end{array}
}{
\begin{array}{c}
(w-i\tau )^\gamma 
\end{array}
} 
\end{array}
& \qquad \qquad \tau >0, 
\end{array}
\end{equation}
where indices $\pm $ denote the cut orientations (0,$\pm \infty ).$ Let us
consider the integrals 
\begin{equation}
\label{A2}I=\lim _{\tau \rightarrow 0+}\int_{-\infty }^{+\infty
}h(a,x-y)h(b,y-z)dy 
\end{equation}
for various combinations of the cut orientations and their locations above
or below the real axis. First let us assume 
\begin{equation}
\label{A2+}a<1,\ b<1,\ a+b>1. 
\end{equation}
For the more general integral (\ref{eq50}) we have shown the integral
vanishes, when the integration line does not separate off the both
singularities, which is the case of the eight integrals (\ref{A2}) involving
combinations $h^{\pm }h_{\pm },\ h^{\pm }h_{\mp },\ h_{\pm }h^{\pm },\
h_{\mp }h^{\pm }$ . Let us calculate (\ref{A2}) when e.g. $x<z$ and the both
singularities are above the real axis. Then 
\begin{equation}
\label{A3}\left| I\right| =\left| \exp (ia\pi )I_1+I_2+\exp (-ib\pi
)I_3\right| =0, 
\end{equation}
where 
\begin{equation}
\label{A3+}I_k=\int_{L_k}\frac{dy}{\left| (x-y)^a(y-z)^b\right| }\qquad
L_1\equiv (-\infty ,x),\ L_2\equiv (x,z),\ L_3\equiv (z,+\infty ). 
\end{equation}
Using simple substitutions in the known relation 
\begin{equation}
\label{A5}\int_0^1x^{\lambda -1}(1-x)^{\mu -1}dx=\frac{\Gamma (\lambda
)\Gamma (\mu )}{\Gamma (\lambda +\mu )} 
\end{equation}
and denoting $d^{a+b-1}\equiv \left| (x-z)^{a+b-1}\right| $ one can get 
$$
I_1=\frac 1{d^{a+b-1}}\cdot \frac{\Gamma (1-a)\Gamma (a+b-1)}{\Gamma (b)} 
$$
\begin{equation}
\label{A4}I_2=\frac 1{d^{a+b-1}}\cdot \frac{\Gamma (1-a)\Gamma (1-b)}{\Gamma
(2-a-b)}\ \ \ \ \ \ 
\end{equation}
$$
I_3=\frac 1{d^{a+b-1}}\cdot \frac{\Gamma (1-b)\Gamma (a+b-1)}{\Gamma (a)}. 
$$
After inserting into (\ref{A3}) we get 
\begin{equation}
\label{A6a}\cos (a\pi )\frac{\Gamma (1-a)\Gamma (a+b-1)}{\Gamma (b)}+\frac{%
\Gamma (1-a)\Gamma (1-b)}{\Gamma (2-a-b)}+\cos (b\pi )\frac{\Gamma
(1-b)\Gamma (a+b-1)}{\Gamma (a)}=0 
\end{equation}
\begin{equation}
\label{A6b}\frac{\sin (a\pi )\Gamma (1-a)}{\Gamma (b)}=\frac{\sin (b\pi
)\Gamma (1-b)}{\Gamma (a)}. 
\end{equation}
The last identity also follows from the known formula (see e.g.\cite{gamma}%
,p.256) 
\begin{equation}
\label{A7}\Gamma (\gamma )\Gamma (1-\gamma )\sin (\gamma \pi )=\pi 
\end{equation}
Now let us calculate (\ref{A2}) for remaining combinations $h^{\pm }h^{\pm
},\ h^{\pm }h^{\mp },\ h_{\pm }h_{\pm },\ h_{\mp }h_{\pm }.$ For $\tau
\rightarrow 0$ it holds 
\begin{equation}
\label{A8}\lim _{\tau \rightarrow 0+}h(\gamma ,w)=\frac{\exp (i\pi c)}{%
\left| w^\gamma \right| }, 
\end{equation}
where the phases $c$ of functions $h$ entering the integral (\ref{A2}) are
in accordance with the convention (\ref{eq6}) summarized in Tab.\ref{tbA1}
This table enables to obtain the phases $c_k$ of their products which are
listed in the first three columns of Tab.\ref{tbA2} 
\begin{table}
\begin{center}
  \begin{tabular}{|c|c|c|} \hline
     {          }  &  {\bf $h(a,x-y)$}  &  {\bf $h(b,x-z)$} \\ 
     {          }  & {\bf $y<x  \qquad   y>x$} & {\bf $y<z  \qquad   y>z$}\\ \hline
$h^{+}$          &  $0 \qquad  \qquad  -a$             &  $-b \qquad \qquad   0$  \\
$h^{-}$          &  $0  \qquad \qquad  -a$              &  $-b \qquad \qquad   0$  \\ \hline
$h_{+}$          &  $-2a \qquad   -a$             &  $-b \qquad   -2b$\\
$h_{-}$           &  $0 \qquad   \qquad  a$               &  $b \qquad \qquad    0$\\ \hline
  \end{tabular}
\end{center}
  \caption{The phases $c$ depending on the cut location\label{tbA1}}
\begin{center}
  \begin{tabular}{|c|c|c|c|c|} \hline
                   &  {\bf $y<x$}  &  {\bf $x<y<z$} & {\bf $z<y$} & {\bf $I$} \\ \hline
$h^{+}h^{+}$&$-b$     &$-a-b$&$-a$&$+\exp (-i\pi [a+b])G$\\
$h_{+}h_{+}$&$-2a-b$&$-a-b$&$-a-2b$&$-\exp (-i\pi [a+b])G$\\
$h^{+}h^{-}$&$-b$      &$-a-b$&$-a$&$+\exp (-i\pi [a+b])G$\\
$h_{+}h_{-}$&$-2a+b$&$-a+b$&$-a$&$-\exp (-i\pi [a-b])G$\\
$h^{-}h^{+}$&$-b$     &$-a-b$&$-a$&$+\exp (-i\pi [a+b])G$\\
$h_{-}h_{+}$&$-b$     &$+a-b$&$+a-2b$&$-\exp (+i\pi [a-b])G$\\
$h^{-}h^{-}$&$-b$      &$-a-b$&$-a$&$+\exp (-i\pi [a+b])G$\\
$h_{-}h_{-}$&$+b$     &$+a+b$&$+a$&$-\exp (+i\pi [a+b])G$\\ \hline
  \end{tabular}
\end{center}
\begin{center}
  \begin{tabular}{|c|c|c|c|c|} \hline
                   &  {\bf $y<z$}  &  {\bf $z<y<x$} & {\bf $x<y$} & {\bf $I$} \\ \hline
$h^{+}h^{+}$&$-b$     &$0$&$-a$&$-G$\\
$h_{+}h_{+}$&$-2a-b$&$-2a-2b$&$-a-2b$&$+\exp (-2i\pi [a+b])G$\\
$h^{+}h^{-}$&$-b$      &$0$&$-a$&$-G$\\
$h_{+}h_{-}$&$-2a+b$&$-2a$&$-a$&$+\exp (-2i\pi a)G$\\
$h^{-}h^{+}$&$-b$     &$0$&$-a$&$-G$\\
$h_{-}h_{+}$&$-b$     &$-2b$&$+a-2b$&$+\exp (-2i\pi b)G$\\
$h^{-}h^{-}$&$-b$      &$0$&$-a$&$-G$\\
$h_{-}h_{-}$&$+b$     &$0$&$+a$&$+G$\\ \hline
  \end{tabular}
\end{center}
  \caption{The  phases  $c_{k}$ of products $h(a,x-y)h(b,y-z)$  and resulting integral $I$ (last column) for the case $x<z$ (upper part) and $x>z$ (lower part) \label{tbA2}}
\end{table}
The integrals of the all combinations summarized in the table can be
expressed like (\ref{A3}), 
\begin{equation}
\label{A9}I=\sum_{k=1}^3\exp (ic_k\pi )I_k 
\end{equation}
If we denote 
\begin{equation}
\label{A10}G\equiv \frac{2i\pi }{d^{a+b-1}}\cdot \frac{\Gamma (a+b-1)}{%
\Gamma (a)\Gamma (b)}, 
\end{equation}
then using identities (\ref{A6a}),(\ref{A7}) the sum (\ref{A9}) can be
evaluated. The results are given in the last column of Tab.\ref{tbA2}. Let
us compare the corresponding rows in upper and lower part of the table.
Obviously for any $x,\ z$ one can write 
\begin{equation}
\label{A14}\lim _{\tau \rightarrow 0+}\int_{-\infty }^{+\infty }h^{\pm
}(a,x-y)h^{\pm }(b,y-z)dy=\lim _{\tau \rightarrow 0+}\frac{2i\pi \Gamma
(a+b-1)}{\Gamma (a)\Gamma (b)}h^{\pm }(a+b-1,x-z) 
\end{equation}
and equally for $h_{\pm }.$ So far we have assumed (\ref{A2+}), however
using identity 
\begin{equation}
\label{A15}\frac d{dw}h(\gamma ,w)=-\gamma h(\gamma +1,w) 
\end{equation}
we can enlarge validity of (\ref{A14}) to any $a,\ b$ 
\begin{equation}
\label{A16}a+b>1,\qquad a,b\neq 0,-1,-2,-3,... 
\end{equation}
In this way we have proven the composition relation (\ref{eq58}). (Note that 
$a=\alpha +1,\ b=\beta +1).$

Alternatively, composition relation can be easily proved in representation
given by Eq. (\ref{eq27}) for $\alpha ,\ \beta <0.$ The relations 
\begin{equation}
\label{A17}D_{\pm }^{\alpha +\beta }(x-z)=\int_{-\infty }^{+\infty }D_{\pm
}^\alpha (x-y)D_{\pm }^\beta (y-z)dy,\qquad \alpha ,\ \beta <0 
\end{equation}
after inserting from (\ref{eq27}) and simple substitutions immediately
follow from (\ref{A5}).
\newpage
{\bf Acknowledgement.~}
I would like to express my gratitude to P. Kol\'a\v r for critical 
reading of the manuscript and valuable comments.

\end{document}